\documentstyle{article}

\makeatletter

\@addtoreset{equation}{section}

\makeatother

\newtheorem{lm}{Lemma}[section]
\newtheorem{pr}{Proposition}[section]
\newtheorem{thm}{Theorem}[section]

\begin{document}

\title{ $L$ operators  and Drinfeld's generators}
\author{ Norifumi Hayaishi and Kei Miki\\
Department of Mathematics,\\
Graduate School of Science, Osaka University\\
Toyonaka 560, Japan}
\date{}
\maketitle

\begin{abstract}
Utilizing the multiplicative formula of universal $R$ matrix, 
the correspondence between the $L$ operators and
 Drinfeld's generators is  explicitly
calculated for quantum group $U_q(g)$ with
 $g=A_l^{(1)}, B_l^{(1)}, C_l^{(1)}, D_l^{(1)}$.
\end{abstract}

\def\hh{h^*}
\def\La{\Lambda}
\def\Q{{\bf Q}}\def\C{{\bf C}}\def\Z{{\bf Z}}
\def\diag{{\rm diag}}\def\End{{\rm End}}
\def\wt{{\rm wt}}\def\ad{{\rm ad}}
\def\tr{{\rm tr}\,}\def\l{{\tilde l}}
\def\A{{\cal A}}\def\R{{\cal R}}\def\L{{\cal L}}
\def\I{{\cal I}}\def\J{{\cal J}}
\def\F{{\cal F}}\def\U{{\cal U}}
\def\T{{\cal T}}\def\D{{\cal D}}\def\N{{\cal N}}\def\R{{\cal R}}
\def\tA{{\tilde \A}}\def\tad{\widetilde{\ad}}\def\th{{\tilde h}}
\def\bg{\bar{g}}\def\bU{\bar{U}}\def\bP{\bar{P}}\def\bQ{\bar{Q}}
\def\bW{\overline W}\def\tW{\widetilde{W}}
\def\hLL{{\hat \L}}\def\hL{{\hat L}}\def\hT{{\hat \T}}
\def\comega{\omega^\vee}
\def\t{{}^t\!}\def\la{\langle}\def\ra{\rangle}\def\n{\nonumber\\}


\section{Introduction}

By Drinfeld  
\cite{Dr1},~\cite{Dr} and Jimbo \cite{J1},
quantum group was invented  by quantizing  
 the  universal enveloping algebras of Kac--Moody algebras. 
 It was defined by deforming the defining relations 
for  Chevally generators. Later Drinfeld  found  
another realization of quantum group
for affine Kac--Moody algebras in \cite{Dr2}.
Finally, Faddeev, Reshetikhin and Takhtatjan \cite{FRT} 
defined quantum group by the $RLL=LLR$ relations,
which deforms the matrix realization of finite dimensional
 simple Lie algebras.  This was later generalized to 
the affine case by Reshetikhin and Semenov-Tian-Shansky \cite{RS}.
Up to now,  the relations among these three definitions are
 clarified as follows. Universal $R$ matrix  by Drinfeld \cite{Dr}
gives the map from the third  formulation to
 the first \cite{FRT}, \cite{RS}.
In the original paper \cite{Dr2},
the isomorphism from the first definition  to the second is given. 
(See \cite{Jing} for the explicit verification.)
Later  in the nontwisted 
affine case,  the map in the opposite
direction is obtained  in \cite{KT2}, without the use of the Braid group,
and in \cite{Be1} by using the Braid group actions.
The correspondence between  the second and  third definitions  is found  
by Ding and Frenkel \cite{DF} for $g=\widehat{gl_n}$.
The purpose of this paper is to generalize the last result
by utilizing the multiplicative formula of universal $R$ matrix.
(See 
\cite{KLP} for a similar calculation in the case $\widehat{DY(sl_2)}$.)
This formula in the nontwisted affine case 
is already obtained in \cite{KT1}
 without the use of Braid group. In this paper, we follow the approach 
by Beck  \cite{Be1}, \cite{Be2} and Lusztig \cite{Lt}. Their results easily 
give the multiplicative formula, and this formula immediately yields 
the desired  correspondence. The  results in this paper would be well known
conceptually. We only hope that the explicit expressions
 here help to clarify the problem.

This paper is organized as follows. 
In section 2, we fix our notations. 
In section 3, we determine the defining relations of quantum group
in the  $L$ operator formalism for $g=A_l^{(1)}$, $B_l^{(1)}$, $C_l^{(1)}$,
 $D_l^{(1)}$, $A_{2l}^{(2)}$, $A_{2l-1}^{(2)}$,  $ D_{l+1}^{(2)}$,
since we cannot find the proof in literature. 
Here we also consider the twisted case since the proof is essentially 
the same as the one in the nontwisted case.
In section 4, after a review of the result by Beck, 
we write down  the multipilcative formula of universal $R$ matrix 
for nontwisted affine Kac--Moody algebras.
In section 5, 
utilizing the result in section 3,  the multiplicative  formula and
the result by Beck,  we calculate the correspondence 
between the $L$ operators and Drinfeld's generators
 for $g=A_l^{(1)}$, $B_l^{(1)}$, $C_l^{(1)}$, $D_l^{(1)}$.
The work by B.  Edwards, announced in \cite{CP}, would immediately 
generalize the results  here to the twisted case.

\section{Notations}

\subsection{} 
Let $g$ be an affine Kac--Moody algebra over $\Q$ and 
$(a_{ij})_{0\le i,j\le l}$   its generalized Cartan matrix.
Let further  $(a_i)_{0\le i\le l}$ and  $(a_i ^\vee)_{0\le i\le l}$ be 
relatively prime positive integers 
such that $ \sum_i a_i^\vee a_{ij}=0$ and 
$\sum_ja_{ij}a_j=0$, respectively. We set $h^\vee=\sum_{i=0}^l a_i^\vee$.
We follow the enumeration  of the vertices of the Dynkin diagram
in  \cite[Chapter 4]{Kac} except in the $A_{2l}^{(2)}$ case, 
in which  we enumerate them  in reverse order.
In our notation,
$a_0=1$ for any $g$ and $a_0^\vee=2$ for $g=A_{2l}^{(2)}$, $=1$ otherwise.
We denote the  finite type  Kac--Moody algebra  corresponding
to the Cartan matrix $(a_{ij})_{1\le i,j\le l}$ by $\bg$.

Let  $\hh$ be a $\Q$ vector space with basis
 $\alpha_0,\cdots,  \alpha_l, \Lambda_0$.
Let further $(d_i)_{0\le i\le l}$ be  relatively prime positive integers
such that $d_i a_{ij}=d_j a_{ji}$.
We introduce a symmetric bilinear form
$(\,\mid \,):\hh\times \hh \rightarrow\Q$ by
$(\alpha_i|\alpha_j)=d_ia_{ij}$, 
$(\alpha_i|\Lambda_0)=d_i\delta_{i0}$ and 
$(\Lambda_0|\Lambda_0)=0$.
Set $\delta=\sum_{i=0}^la_i\alpha_i$ and 
 $\alpha_i^\vee=\alpha_i/d_i$ $(0\le i\le l)$.
Define $\omega_i$ and  $\omega_i^\vee \in \hh$  $(1\le i\le l)$  by
$(\alpha_j^\vee|\omega_i)=\delta_{ij} \quad(1\le j\le l)$,
 $(\delta|\omega_i)=(\Lambda_0|\omega_i)=0$
 and $\omega_i^\vee=\omega_i/d_i$.
We further define several $\Z$ submodules  of $\hh$ by
$Q=\oplus_{i=0}^l\Z \alpha_i$, $\Sigma=Q\oplus\Z\Lambda_0/d_0$, 
$\Gamma=\Sigma+p\Z\omega_1$ and  $\bar{Q}=\oplus_{i=1}^l\Z\alpha_i$.
Here and hereafter $p=2$ for $g=A_2^{(2)}$ and $=1$ otherwise.
Finally, as usual, we introduce an order relation on $Q$ by
$\lambda\ge\mu \Leftrightarrow \lambda-\mu\in Q_+ 
=\oplus_{i=0}^l\Z_{\ge 0}\alpha_i$.

\subsection{} 
Let $q$ be an  indeterminate and set  $F=\Q(q)$.
Following \cite{Dr1},~\cite{J1}, 
let $\U=U^\Gamma_q(g)$ be the  $F$ algebra generated by
$e_i$, $f_i$ $(0\le i\le l)$ and  $k_\lambda$ $(\lambda\in \Gamma)$
with the relations,

\begin{eqnarray}
&&k_\lambda k_\mu=k_{\lambda+\mu},\quad k_0=1,\\
&&k_\lambda e_i k_{-\lambda}=q^{(\lambda|\alpha_i)} e_i,\quad 
k_\lambda f_i k_{-\lambda}=q^{-(\lambda|\alpha_i)} f_i, \\
&&[e_i,f_j]=\delta_{ij}{k_i-k_i^{-1}\over q_i-q_i^{-1}},\\
&&\sum_{r=0}^{1-a_{ij}} (-1)^r
e_i^{(r)}e_je_i^{(1-a_{ij}-r)}=0, \quad (i\ne j)\\
&&\sum_{r=0}^{1-a_{ij}} (-1)^r
f_i^{(r)}f_j f_i^{(1-a_{ij}-r)}=0,\quad (i\ne j)
\end{eqnarray}
where  $k_i^{\pm 1}=k_{\pm \alpha_i}$ and 
$$
q_i=q^{d_i},\quad 
x_i^{(n)}={x_i\over [n]_i!}\quad (x=e,f), \quad 
[n]_i={q_i^n-q_i^{-n}\over q_i-q_i^{-1}},\quad [n]_i!=\prod_{k=1}^n[k]_i.
$$
This algebra has the following Hopf algebra structure,
\begin{eqnarray}
&&
\Delta(e_i)=e_i\otimes 1+k_i\otimes e_i,\quad 
\Delta(f_i)=f_i\otimes k^{-1}_{i}+ 1 \otimes f_i,\quad
\Delta(k_\lambda)=k_\lambda\otimes k_\lambda, \n
&&
S(e_i)=-k^{-1}_i e_i,\quad 
S(f_i)=-f_ik_i,\quad 
S(k_\lambda)=k_{-\lambda}, \n
&&
\epsilon(e_i)=\epsilon(f_i)=0,\quad \epsilon(k_\lambda)=1.
\end{eqnarray}
We define several subalgebras of $\U$ by
\begin{eqnarray*}
&&\U^+=\la e_i|0\le i\le l\ra, \quad \U^-=\la f_i|0\le i\le l\ra,
\quad \U^0=\la k_\lambda | \lambda\in \Gamma \ra,\\
&&\U^{\ge 0}=\U^0 \U^+,\quad \U^{\le 0}=\U^{-}\U^0.
\end{eqnarray*}
Then $\U^\pm=\oplus_{\mu\in Q_+} \U^\pm_{\pm \mu}$ where
$$
\U^\pm_\mu=\{x\in \U^\pm |\,k_\lambda x k_{-\lambda}
=q^{(\lambda|\mu)}x\quad (\forall \lambda\in \Gamma)\}
\quad (\mu\in Q).
$$
Moreover the following is  known.
\begin{pr}\label{pr:triang} \cite{Ro}, \cite{Ls}, \cite{Lt}

(1) The multiplication map  $\U^-\otimes \U^0\otimes \U^+\to  \U$
induces an isomorphism of  vector spaces.

(2) $\{k_\lambda|\lambda\in \Gamma\}$  forms a basis of $\U^0$.

(3) $\dim_F \U^\pm_{\pm \mu}={\cal P}(\mu)$ for $\mu \in Q$,
where ${\cal P}(\mu)$ is Kostant's partition  function for $g$.
\end{pr}
Similarly we define $U=U_q^\Sigma(g)$, 
$U'=U_q^{Q}(g)$ and $\bU=U_q^{\bar{Q}}(\bg)$. 
We   identify $U$ (resp. $U'$)   with the subalgebra of $\U$ (resp. $U$).
In particular, $\U^\pm=U^\pm=U'{}^\pm$.
Finally, for a $\bU$ module $M$ and 
$\mu\in \bar{P}=\oplus_{1\le i\le l}\Z \omega_i$,
we define its weight space $M_\mu$ by
$$
M_\mu=\{v\in M\mid k_\lambda v=q^{(\lambda|\mu)}v
\quad (\forall \lambda \in \bar{Q})\}.
$$

\subsection{}
Hereafter to the end of the next section, we consider
$g=A_l^{(1)}$ $(l\ge 1)$, $B_l^{(1)}$ $(l\ge 3)$ ,
 $C_l^{(1)}$ $(l\ge 2)$, $D_l^{(1)}$ $(l\ge 4)$, $A_{2l}^{(2)}$ $(l\ge 1)$,
 $A_{2l-1}^{(2)}$ $(l\ge 3)$, $ D_{l+1}^{(2)}$ $(l\ge 2)$.
Let $V(p \omega_1)$ be the  irreducible highest weight module  with highest 
weight $p \omega_1$ of $\bU$. Let further $v_i$ 
$(1\le i\le \dim V(p \omega_1))$ be 
a basis of  this module consisting of weight vectors
 such that  $i<j$ if $\eta_i>\eta_j$,
where $\eta_i\in \bQ$ is the weight of $v_i$ as a $\bU$ module.
$V(\omega_1)\oplus F$ in the case $g= D_{l+1}^{(2)}$,
 and $V(p \omega_1)$  otherwise can be made into  representations
  of $U'$, on which $k_{\pm \delta}$ acts as $1$.
We shall denote them by $(\rho, V)$ and set $N=\dim V$.
In the case $g=D_{l+1}^{(2)}$,   by $v_N$
we denote   the basis vector of the trivial representation  as a $\bU$ module.
The representations $(\rho, V)$ are  
the same as the ones used in \cite{J2} and,
in the limit $q=1$, reduce to the ones in \cite[Chapters 7 and 8]{Kac}. 
Let $\tau_a$ $(a\in F^{\times})$ be the automorphism of $\U$ determined by
$\tau_a(e_i)=a^{\delta_{i0}}e_i$, $\tau_a(f_i)=a^{-\delta_{i0}}f_i$ and  
$\tau_a(k_\lambda)=k_\lambda$.
Put $\rho_a=\rho\circ\tau_a|_{U'}$. 
Set $M=l+1$ and $a_i=q^{-2(i-1)}$$(1\le i\le M)$ for $g=A_l^{(1)}$. 
In other cases, set 
$M=2$ and $a_1=1$, $a_2=\sigma q_0^{-h^\vee\!/a_0^\vee}$.
Here $\sigma=-1$ for $g=A_{2l}^{(2)}$, $A_{2l-1}^{(2)}$ and $=1$ otherwise.
Then in $(\rho_{a_1}\otimes \cdots \rho_{a_M}, V^{\otimes M})$, 
there exists a trivial representation of $U'$.
 We denote its basis vector by $w$.
In the case $M=2$,  we define a matrix $J=(J_{ij})_{1\le i,j\le N}$ by 
$w=\sum_{i,j}J_{ij} v_i\otimes v_j$.
See Appendix for the  explicit expressions of  $N$,
 $\eta_i$, $(\rho, V)$ and $w$. 
Below we identify $\End V$ with $M_N(F)$  through 
the basis $v_1,\cdots, v_N$ normalized as in Appendix.

\subsection{}
Let  $\Theta$ be the quasi-universal $R$  matrix 
 of $U$ \cite{Dr},~\cite{T},~\cite{Lt}.
This is an invertible element of (an appropriate completion of )
 $U\otimes U$ uniquely determined by
\begin{eqnarray}
&&\Theta=\sum_{\mu\ge0} \Theta_\mu,\quad 
 \Theta_\mu\in U^+_\mu\otimes U^-_{-\mu},\n
&&\Theta_0=1\otimes 1,
\quad \sigma\circ \Delta(x)\Theta
=\Theta \Psi\circ\Delta(x)\quad (x\in U)\label{eq:r1}
\end{eqnarray}
Here $\sigma(x\otimes y)=y\otimes x$, and 
$\Psi$ is the  automorphism of $U\otimes U$ defined by
\begin{eqnarray*}
&&\Psi(e_i\otimes 1)=e_i\otimes k^{-1}_i,\quad 
\Psi(f_i\otimes 1)=f_i\otimes k_i,\n
&&\Psi(k_\lambda\otimes k_\mu)
=k_\lambda\otimes k_\mu,\quad \Psi=\sigma\circ \Psi
\end{eqnarray*}
This element further satisfies,
\begin{eqnarray}
&&\Theta_{\alpha_i}=-(q_i-q_i^{-1})e_{i}\otimes f_{i}\quad (0\le i\le l),
\label{eq:r2}\\
&&(1\otimes \Delta)\Theta=\Theta_{13}\Psi_{13}(\Theta_{12}),\quad
(\Delta\otimes 1)\Theta=\Theta_{13}\Psi_{13}(\Theta_{23}),\label{eq:r3}\\
&&\Theta_{12}\Psi_{12}(\Theta_{13})\Theta_{23}=
\Theta_{23}\Psi_{23}(\Theta_{13})\Theta_{12},
\end{eqnarray}
where $\Theta_{13}$, for example, 
denotes $\sum_{i} a_i\otimes 1\otimes b_i$ 
for $\Theta=\sum_i a_i\otimes b_i$
as usual.

For $\Theta^{-1}$, define  $(\Theta^{-1})_\mu$ 
$(\mu\in Q_+)$ as for $\Theta$. Let us introduce  
   $\hLL^\pm(z)\in (M_N(F) \otimes U^{\mp})[[z^{\pm 1}]]$
and  $\R(z) \in (M_N(F)\otimes M_N(F))[[z]]$ by
\begin{eqnarray}
&&\hLL^+(z)=\sum_{\mu\ge 0} (\rho \otimes 1)
(\Theta_\mu) z^{{(\Lambda_0|\mu)\over d_0}},\quad 
\hLL^-(z)=\sum_{\mu\ge 0} (\rho \otimes 1)\left(
((\Theta_{21})^{-1})_\mu\right) z^{-{(\Lambda_0|\mu)\over d_0}}, \n
&&\R(z)=(1\otimes \rho)\L^+(z).
\end{eqnarray}
We further introduce  $L$ operators and $R$ matrix by
\begin{eqnarray}
&&\hL^+(z)=\hLL^+(z)\hT^{-1},\quad \hL^-(z)=\hT \hLL^{-}(z),\n
&& R(z)=\R(z) T^{-1}=\sum_{n\ge 0}
\sum_{1\le i_1, i_2, j_1, j_2 \le N} R[n]^{i_1 j_1}_{i_2 j_2}
E_{i_1 i_2}\otimes E_{j_1 j_2} z^n,
\end{eqnarray}
where 
$$
\hT^{\pm 1}=\sum_{1\le i\le N} E_{ii}\otimes k_{\pm \eta_i},\quad
T^{\pm 1}=\sum_{1\le i, j\le N}  
q^{\pm( (\eta_i|\eta_j)-(\eta_1|\eta_1))}E_{ii}\otimes E_{jj},
$$
and $E_{ij}$ is a matrix unit.

Finally we state the simple properties of $R$ matrix, 
which we shall need in the next section. 
Let $(\dot{\rho}, \Q^N)$ denote  the representation of $g'=[g,g]$ obtained
by taking the limit $q=1$ of $\rho(\cdot)\in M_N(F)$.
Let $r(z)\in (\dot{\rho}(g')\otimes_\Q\! \dot{\rho}(g'))[[z]]$
be the trigonometric classical $r$ matrix \cite{BD} of type $g$ normalized 
by the condition that the  diagonal part of the $z^0$ component of $r(z)$ is 
$\sum_{1\le i,j\le N}{(\eta_i|\eta_j)\over 2} E_{ii}\otimes E_{jj}(=:r_0)$.
Let  further $A$ be the  $\Q$ subalgebra of $F$ consisting 
of the elements which have no pole at $q=1$.
\begin{lm}
\label{lm:R}
\begin{eqnarray*}
&&(1)\, R[0]^{i_1 j_1}_{i_2 j_2}\ne 0 
\mbox{  only if  } i_1 \le i_2 \mbox{  and  } j_1 \ge j_2\\
&&(2)\, \eta_{i_1}+\eta_{j_1}=\eta_{i_2}+\eta_{j_2},
\quad  i_1=i_2 \Leftrightarrow j_1=j_2 
\mbox{   when   } R[n]^{i_1 j_1}_{i_2 j_2}\ne 0 \\
&&(3)\,  \R(z)=1-(q-q^{-1})(r(z)-r_0)
\quad  {\rm mod} \quad (q-q^{-1})^2 (M_N(A)\otimes_A \! M_N(A))[[z]]
\end{eqnarray*}
\end{lm}
Proof.
(3) Let $r_\mu$ $(\mu>0)$ be the elements
 of $M_N(\Q)\otimes_\Q\! M_N(\Q)$
uniquely determined by the following conditions,
\begin{eqnarray*}
&& r_{\alpha_i}=d_i\dot{\rho}(e_i)\otimes \dot{\rho}(f_i)
 \quad (0\le i\le l),\\
&&[1\otimes \dot{\rho}(e_i), r_\mu ]+
[\dot{\rho}(e_i)\otimes 1, r_{\mu-\alpha_i}]=0
 \quad (\mu\ne \alpha_i>0, 0\le i\le l),\\
&&[\dot{\rho}(f_i)\otimes 1, r_\mu]+
[1\otimes \dot{\rho}(f_i), r_{\mu-\alpha_i}]=0 
\quad (\mu\ne \alpha_i>0, 0\le i\le l),\\
&&\sum_{k=1}^M(r_\mu)_{k\,  M+1}
\dot{w}\otimes 1_{\scriptscriptstyle {\Q^N}} =0\quad (\mu>0)
\quad \mbox{for } g\ne A_{2l}^{(2)}, A_{2l-1}^{(2)},\\
&&({\rm tr}\otimes 1) (r_{\mu})=0 
\quad (\mu>0)\quad \mbox{for } g= A_{2l}^{(2)}, A_{2l-1}^{(2)},
\end{eqnarray*}
where  $\dot{w}=w|_{q=1}$.
Then $r(z)=\sum_{\mu\ge 0} r_\mu z^{(\Lambda_0|\mu)\over d_0}$.
Hence  by considering  the equations from (\ref{eq:r1}), (\ref{eq:r2}) 
and  (\ref{eq:r3}) in the limit $q=1$,
the claim  can be easily proven by induction on ${\rm ht}\mu$.
 In the case $g=A_{2l}^{(2)}$, $A_{2l-1}^{(2)}$, 
we also use (\ref{eq:trace}) below with 
$l_{ij}[-n]=-\hLL^+_{ij}[-n]/(q-q^{-1})$.
(In the nontwisted case,  it is also possible to 
use the multiplicative formula of $\Theta$ below.
See \cite[Section 6]{Be1} and \cite[Section 1]{Lss}.)


\section{$L$ operator formalism}

\subsection{}
Let $g$, $(\rho, V)$, $N$, $\eta_i$ $(1\le i\le N)$, $R(z)$, 
$M$, $a_i$ $(1\le i\le M)$ and  $w$  as before.
Following \cite{FRT},~\cite{RS},~\cite{DF},~\cite{FR}, let 
$\A=\A(g)$ be the  $F$ algebra generated by $L^+_{ij}[0]$
 $(1\le i\le j\le N)$, $L^-_{ij}[0]$ $(1\le j\le i\le N)$,
$L^\pm_{ij}[\mp n]$ $(1\le i,j \le N,n\in \Z_{>0})$,
 $C^{\pm 1}$ and $D^{\pm1}$ with the relations
\begin{eqnarray}
&&L^\pm_{ii}[0]L^\mp_{ii}[0]
=C^{\pm 1}C^{\mp1}=D^{\pm 1}D^{\mp 1}=1\label{eq:a1}\\
&&L^\pm_{ii}[0]=1\mbox{  for $i$ such that $\eta_i=0$}\label{eq:a2}\\
&&C^{\pm 1} \mbox{ central}\\
&&DL^\pm_{ij}[\mp n]D^{-1}=q^{\mp n}L^{\pm}_{ij}[\mp n]\\
&&R_{12}(Cz/w)L_1^+(z)L_2^-(w)=L_2^-(w)L_1^+(z)R_{12}(C^{-1}z/w)\\
&&R_{12}(z/w)L_1^\pm(z)L_2^\pm(w)=L_2^\pm(w)L_1^\pm(z)R_{12}(z/w)\\
&&L^\pm_1(a_1 z)\cdots L^\pm_M(a_M z) w\otimes 1=w\otimes 1\label{eq:lw}\\
&&(G\otimes 1)L^\pm(z) (G^{-1}\otimes 1)=L^\pm(-z) \quad \mbox{for }
g=D_{l+1}^{(2)} \label{eq:lD}
\end{eqnarray}
Here
$L^\pm(z)$ is defined in terms of the generating functions 
$L^\pm_{ij}(z)=\sum_{n\ge 0}L^\pm_{ij}[\mp n] z^{\pm n}$ as 
$$
L^\pm(z)=\sum_{1\le i,j\le N} E_{ij}\otimes L^\pm_{ij}(z)
 \in (M_N(F)\otimes \A)[[z^{\pm1}]], 
$$
$L_i^\pm(z)$ denotes $L^\pm(z)$ whose matrix part acts
 on the $i$ th space,  and
$$
G=\diag(1,\cdots,1,-1)\in M_N(F) \quad  \mbox{for   } g=D_{l+1}^{(2)}.
$$
It is well known and easy to check that $\A$ has 
the following Hopf algebra structure
\begin{eqnarray}
&&\Delta(L^{+}_{ij}(z))
=\sum_{k=1}^N L^{+}_{kj}(C_2^{-1}z)\otimes L^{+}_{ik}(z),\quad
\Delta(L^{-}_{ij}(z))=\sum_{k=1}^N L^{-}_{kj}(z)
\otimes L^{-}_{ik}(C_1^{-1}z),\n
&&\Delta(C^{\pm 1})=C^{\pm 1}\otimes C^{\pm 1},
\quad \Delta(D^{\pm 1})=D^{\pm 1}\otimes D^{\pm 1},\n
&&\epsilon(L^\pm_{ij}(z))=\delta_{ij},\quad  
\epsilon(C^{\pm 1})=\epsilon(D^{\pm 1})=1,\n
&&S(L^{\pm}(z))=\t(( \t L^\pm(Cz))^{-1}),\quad S(C^{\pm 1})=C^{\mp 1},
\quad S(D^{\pm 1})=D^{\mp 1},\label{eq:copr}
\end{eqnarray}
where $C_1=C\otimes 1$,  $C_2=1\otimes C$ and ${}^t$ denotes 
the matrix transpose in the first space. 
It is also well known and easy to check that
there is a surjective Hopf algebra homomorphism
 $\phi:\A\to \U$ defined by
$$
L^\pm(z)\mapsto\hL^\pm(z), \quad C^{\pm 1}\mapsto k_{\pm \delta},
\quad D^{\pm 1}\mapsto k_{\pm \Lambda_0/d_0}.
$$
In the next subsection, we shall show the injectivity of $\phi$ 
by the standard specialization argument at $q=1$ \cite{Ls},~\cite{Lss}. 
The only subtlety appears in the case $A_{2l}^{(2)}$, $A_{2l-1}^{(2)}$.

\subsection{}

Set 
\begin{equation}
\L^+(z)=L^+(z)\T,\quad \L^-(z)=\T^{-1} L^{-}(z),
\quad \T^{\pm 1}=\sum_{1\le i\le N} E_{ii}\otimes L^{\mp}_{ii}[0],
\end{equation}
and define the components $\L^{\pm}_{ij}[\mp n]$ as before.
We further define  subalgebras of $\A$ by
$$
\A^{\pm}=\la\L^{\mp}_{ij}[\mp n] |1\le i,j\le N, n\ge 0\ra, \quad 
\A^0=\la L^{\pm}_{ii}[0], C^{\pm 1}, D^{\pm 1} |1\le i\le N\ra, 
$$
and,  for $\mu\in Q$,  we   set
$$
\A^\pm_{\mu}=\{x\in \A^\pm | L^{-}_{ii}[0] x L^{+}_{ii}[0]
=q^{(\eta_i|\mu)} x\quad (1\le i\le N),\quad 
D x D^{-1}=q^{{(\Lambda_0|\mu)\over d_0}} x \}.
$$

\begin{lm}\label{lm:A}

(1) $L^{-}_{kk}[0] L^\pm_{ij}[\mp n]
=q^{-(\eta_k|\eta_i-\eta_j)}  L^\pm_{ij}[\mp n] L^{-}_{kk}[0]$.

(2) $\A^0 $ is  commutative.

(3) $\prod_{i=1}^{N}{L^\pm_{i i}[0]}^{n_i}=1
 \mbox{   if    }\sum_i n_i \eta_{i}=0$.

(4) $\A=\A^-\A^0\A^+$.

\end{lm}

Proof. (3) follows from (\ref{eq:a2}) 
and the $z^0$ components  of  (\ref{eq:lw}) and (\ref{eq:lD}).

\begin{lm}\label{lm:dim}
$\A^{\pm}=\oplus_{\mu \in Q_+} \A^\pm_{\pm \mu}$
and $\dim_F \A^\pm_{\pm \mu}\le {\cal P}(\mu)$.
\end{lm}

Proof. 
This can be   proven  as in \cite[Proposition 1.13]{Lss}.

Let $\N$ be the $F$ algebra generated by $\L_{ij}[0]$
$(1\le i< j\le N)$ and $\L_{ij}[-n]$ $(1\le i,j\le N, n>0)$
with the relations
\begin{eqnarray}
&&\R_{12}(z/w)T_{12}^{-1}\L_1(z) T_{12}\L_2(w)
=\L_2(w)T_{12}\L_1(z)T_{12}^{-1}\R_{12}(z/w)\label{eq:n1}\\
&&\L_1(a_1 z)T_{12}\L_2(a_2z)T_{12}^{-1}\cdots 
T_{1;M}\L_M(a_Mz)T_{1;M}^{-1}w\otimes 1=w\otimes 1
\label{eq:n2}\\
&&(G\otimes 1)\L(z)(G^{-1}\otimes 1)=\L(-z)\quad
\mbox{for   } g=D_{l+1}^{(2)}\label{eq:n3}
\end{eqnarray}
where $T_{1;k}=T_{1k}\cdots T_{k-1\,k}$,
and  $\L(z)$  are defined as $L^+(z)$  with $\L_{ii}[0]=1$.
Next let  $\N_A$ be the $A$ algebra generated by 
$l_{ij}[0]$ $(1\le i< j\le N)$ and 
$l_{ij}[-n]$ $(1\le i,j\le N, n>0)$ with the following relations.
Let $l(z)$ be defined as before with  $l_{ii}[0]=0$. 
Firstly substitute  $\L(z)=1-(q-q^{-1})l(z)$ 
into (\ref{eq:n1}--\ref{eq:n3}) and 
then  divide (\ref{eq:n1}) by $(q-q^{-1})^2$ and the rest  by $q-q^{-1}$.
We take the thus obtained  equations as the defining relations of $\N_A$.
In the case $g=A_{2l}^{(2)}, A_{2l-1}^{(2)}$, we impose another relation,
\begin{eqnarray}
{\tr l(a_2z)-\tr l(a_2^{-1}z)\over
 q-q^{-1}}&=&\tr\l(z)^2-\tr l(z)^2+ \nonumber\\
(q-q^{-1})&\times &\left(\tr l(z)^2\l(a_2z)-\tr l(a_2^{-1}z)\l(z)^2\right) 
\label{eq:trace}
\end{eqnarray}
where  $\l(z)=\sum_{n\ge 0} \sum_{ij} E_{ij}\otimes
{u_i\over u_j} q^{(\eta_i|\eta_i-\eta_j)}l_{N+1-j\, N+1-i}[-n] z^n$
(see  Appendix A.3   for the $u_i$)
 and ${\rm tr}$ denotes the trace of the matrix part.
Finally let $n^-$ be the subalgebra of $g$ generated by the $f_i$ 
in the notations of \cite{Kac} and $U(n^-)$ its enveloping algebra.
$U(n^-)$ can be defined by the generators $l_{ij}[0]$ 
$(1\le i< j\le N)$, $l_{ij}[-n]$ $(1\le i,j\le N, n>0)$
and the relations
\begin{eqnarray}
&&[r_{12}(z/w), l_1(z)+l_2(w)]+[l_1(z),l_2(w)]+
  [(r_0)_{12}, l_1(z)-l_2(w)]=0,\n
&&l(z)(\dot{J}\otimes 1) +(\dot{J}\otimes 1) \,{}^t l(\sigma z)=0\quad
 \mbox{ for   } g\ne A_l^{(1)},\n
&&\tr l(z)=0 \quad \mbox{ for   } 
g= A_l^{(1)}, A_{2l}^{(2)}, A_{2l-1}^{(2)},\n
&&(G\otimes 1) l(z)(G^{-1}\otimes 1)=l(-z)\quad 
\mbox{ for   } g=D_{l+1}^{(2)},  \label{eq:un}
\end{eqnarray}
where $\dot{J}=J|_{q=1}$. (See \cite[Chapters 7 and 8]{Kac}.)

Set $\N_F=\N_A\otimes_A F$ and $\N_\Q=\N_A\otimes_A\Q$,
where $A$ acts on $F$ naturally and on $\Q$ via $q\mapsto  1$. 
In the case $g= A_{2l}^{(2)}, A_{2l-1}^{(2)}$, 
(\ref{eq:trace}) follows from (\ref{eq:n2}) in $\N$.
Thanks to Lemma \ref{lm:R} (3), in the limit $q=1$
 the relations of $\N_A$ reduce to
 (\ref{eq:un}), cf. \cite{RS},~\cite{FR}.
Therefore 
\begin{eqnarray}
&&\N_\F\simeq\N \quad  (l_{ij}[-n]
\otimes 1 \mapsto -\L_{ij}[-n]/(q-q^{-1})),\n
&&\N_\Q\simeq U(n^{-}) \quad (l_{ij}[-n]
\otimes 1 \mapsto l_{ij}[-n]). \label{eq:NN}
\end{eqnarray}

As easily checked, $\N_A$  is  a  $Q$--graded algebra by assigning
 $\eta_j-\eta_i-n\delta\in Q$ to $l_{ij}[-n]$ and 
each  homogeneous subspace $(\N_A)_{\mu}$
 $(\mu\in Q)$ is finitely generated over $A$.
Moreover  there is a surjective  
$F$ algebra homomorphism from $\N$ to $\A^{-}$
which maps $\L_{ij}[-n]$ to $\L^+_{ij}[-n]$. 
Therefore, thanks to (\ref{eq:NN})
we obtain $\dim_F \A^-_{-\mu}\le \dim_F\, (\N_A)_{-\mu}\otimes_A  \!F\le 
\dim_\Q \,(\N_A)_{-\mu}\otimes_A \!\Q={\cal P}(\mu)$.

\begin{thm}
For 
$g=A_l^{(1)}$, $B_l^{(1)}$, $C_l^{(1)}$, $D_l^{(1)}$, $A_{2l}^{(2)}$,
$A_{2l-1}^{(2)}$,  $ D_{l+1}^{(2)}$,
$\phi:\A(g)\to U_q^\Gamma(g)$ is a Hopf algebra isomorphism.
\end{thm}

Proof. 
Clearly  $\phi(\A^0)=\U^0$ and 
$\phi(\A^\pm_{\pm \mu})=\U^\pm_{\pm \mu}$. Since
$\Gamma=\oplus_{i=1}^l\Z \eta_i\oplus\Z\delta \oplus\Z\Lambda_0/d_0$, 
Proposition \ref{pr:triang} (2), 
 Eq. $\!$(\ref{eq:a1}) and  Lemma \ref{lm:A} (2) (3) 
show that  $\phi|_{\A^0}$ is injective.
Proposition \ref{pr:triang} (3) and Lemma \ref{lm:dim} 
 prove that $\phi|_{\A^\pm}:\A^\pm\to \U^\pm$ are
isomorphisms.
Hence   Proposition \ref{pr:triang}
 (1) and Lemma \ref{lm:A} (4) yield the claim.

\noindent {\it Remark.} In the case $g=A_l^{(1)}$, 
we can replace (\ref{eq:lw}) in the definition of $\A(g)$ by 
the following well known conditions for quantum determinants,
$$
\sum_{\gamma\in {\cal S}_{l+1}} (-q)^{l(\gamma)}
L^\pm_{1\,\gamma(1)}(a_1 z)\cdots L^\pm_{l+1\,\gamma(l+1)}(a_{l+1}z)=1
$$
where $l(\gamma)$ denotes the length function of ${\cal S}_{l+1}$.
For   the above relation is contained in (\ref{eq:lw}) 
and is sufficient to  give Lemma \ref{lm:A} (3) and 
$\tr l(z) =0$ in (\ref{eq:un}).


\section{ Multiplicative formula of $\Theta$}

\subsection{}
In this  section, we shall consider nontwisted 
affine Kac--Moody algebra $g$.
In this subsection, for completeness and a later purpose, 
we review the result by Beck \cite{Be1} 
on the  isomorphism of $U$ and Drinfeld's
 realization $\D$ (see below). We slightly change the notations
in order to do without the square root of the central element.

Let  $\bP^\vee=\oplus_{i=1}^l\Z\omega_i^\vee$ be 
the  coweight lattice of $\bg$ and 
$\bQ^\vee=\oplus_{i=1}^l\Z\alpha_i^\vee$
the coroot lattice.
Let $\bW$ and $W$ be the Weyl groups of $\bg$ and $g$, respectively.  
Set $\tW=\bW|\!\!\!\times  \bP^\vee$, 
where $\bW$ acts on $\bP^\vee$ naturally.
$\tW$ is  known to be isomorphic to $\Pi |\!\!\times W$.
 Here $\Pi$ is the  subgroup of the Dynkin diagram automorphism
 isomorphic to $\bP^\vee/\bQ^\vee$,  and
$\pi\in \Pi$ acts on $s_i$, a simple reflection with respect 
to $\alpha_i$, as $\pi s_i\pi^{-1}=s_{\pi(i)}$.
$\tW$ acts on $Q$. In particular,
$x\in \bP^\vee$  acts  on  $\alpha\in Q$ as 
$x(\alpha)=\alpha-(\alpha|x)\delta$.

Let $T_i$ $(0\le i\le l)$ \cite{Lt} and 
$T_\pi$ $(\pi\in \Pi)$ be the automorphisms  of $U'$ defined by
\begin{eqnarray*}
&&T_i(e_i)=-f_i k_i,\quad T_i(f_i)=-k_i^{-1} e_i,\quad
T_i(k_\lambda)=k_{s_i(\lambda)}\\
&&T_i(e_j)=ad(e_i^{(-a_{ij})})\cdot e_j,
\quad T_i(f_j)=f_j\cdot \tad(f_i^{(-a_{ij})})\quad (i\ne j)\\
&&T_\pi(e_i)=e_{\pi(i)},\quad T_{\pi}(f_i)=f_{\pi(i)},
\quad T_{\pi}(k_{\alpha_i})=k_{\alpha_{\pi(i)}}
\end{eqnarray*}
where  $\ad$ (resp. $\tad$) is a left (resp. right) 
adjoint action  defined by
$$
\ad(x)\cdot y=\sum x^{(1)} y S(x^{(2)}),
\quad y\cdot\tad(x)=\sum S(x^{(1)}) y x^{(2)}
$$
for $x,y \in U$ and $\Delta(x)=\sum x^{(1)}\otimes x^{(2)}$.
For ${\tilde w}\in \tW$, we say an expression 
${\tilde w}=\pi s_{i_1}\cdots s_{i_n}$ $(\pi\in \Pi)$
is reduced if  $s_{i_1}\cdots s_{i_n}$ is so in $W$.
For any reduced expression 
${\tilde w}=\pi s_{i_1}\cdots s_{i_n}$, we can set
$T_{\tilde w}=T_\pi T_{i_1}\cdots T_{i_n}$. Then $T_{{\tilde w}}$ 
$({\tilde w}\in \tW)$ gives a Braid group  action on $U'$.
Let  $\Omega$ be  the $\Q$ algebra 
anti--automorphism of $U'$ defined by
$$
e_i\mapsto f_i,\quad f_i\mapsto e_i,
\quad k_\lambda\mapsto k_{-\lambda},
\quad q\mapsto q^{-1}.
$$
This satisfies $\Omega T_{\tilde w}=T_{\tilde w}\Omega$.

Fix  $o(i)\in \{-1,1\}$ $(1\le i\le l)$ 
 such that $o(i)=-o(j)$ if $a_{ij}<0$.
Let us define $x^{(\pm)}_{i,k}$  and 
$h_{i,r}\in U$ $(1\le i\le l, k\in \Z, r\in \Z-\{0\})$ by
\begin{eqnarray}
&&x^{(+)}_{i,k}=o(i)^k T_{\omega_i^\vee}^{-k}(e_i),
\quad x^{(-)}_{i,k}=\Omega(x^{(+)}_{i,-k}), \\
&&k^{\mp 1}_i\exp\left(\mp(q_i-q_i^{-1})
\sum_{r>0} h_{i,\mp r} z^{\pm  r}\right)
=\sum_{r\ge 0}\phi^\pm _{i,\mp r}z^{\pm r}\label{eq:phi}
\end{eqnarray}
where
\begin{eqnarray}
&& \phi^-_{i,0}=k_i,\quad \phi^-_{i,r}
=(q_i-q_i^{-1})k_{r\delta}[e_i, x^{(-)}_{i,r}] \quad (r>0)\n
&&\phi^+_{i,-r}=\Omega(\phi^-_{i,r})\quad (r\ge 0)
\end{eqnarray}

Let $\D=\D(g)$ \cite{Dr2} be the $F$ algebra generated by 
$x^{(\pm)}_{i,k}$, $h_{i,r}$, $k_i^{\pm 1}$, $C^{\pm 1}$, $D^{\pm 1}$ 
$(1\le i\le l, k\in \Z, r\in \Z-\{0\})$ with the relations,
\begin{eqnarray}
&&k_i^{\pm 1}k_i^{\mp1} =C^{\pm 1}C^{\mp 1}=D^{\pm 1}D^{\mp 1}=1\\
&&C^{\pm 1} \mbox{   central}\\
&&[ k_i, k_j]=[k_i,  D]=[k_i, h_{j,r}]=0, \\
&&[h_{i,r}, h_{j,s}]=\delta_{r+s,0} {[r a_{ij}]_i\over r}
{C^r-C^{-r}\over q_j-q_j^{-1}},\label{eq:h}\\
&& k_i x_{j, k}^{(\pm)} k_i^{-1}
=q^{\pm (\alpha_i|\alpha_j)} x_{j, k}^{(\pm)}, \\
&& D x_{i,k}^{(\pm)} D^{-1}=q^k x_{i,k}^{(\pm) }, 
\quad D h_{i,r} D^{-1}=q^r h_{i, r}, \\
&& [h_{i, r}, x_{j, k}^{(\pm)}]=\pm {[r a_{ij}]_i\over r}
C^{r\mp |r|\over 2} x_{j, r+k}^{(\pm)},\\
&& [x_{i, m}^{(+)}, x_{j, n}^{(-)}]={\delta_{ij}\over q_i-q_i^{-1}}
\left( C^{-n} \phi^{-}_{i, m+n}-C^{-m}\phi^+_{i, m+n} \right), \\
&& x_{i, m+1}^{(\pm)}x_{j,n}^{(\pm)}-q^{\pm (\alpha_i|\alpha_j)}
 x_{j,n}^{(\pm)}x_{i, m+1}^{(\pm)}\n
&&= q^{\pm (\alpha_i|\alpha_j)} x_{i, m}^{(\pm)}
 x_{j, n+1}^{(\pm)}-x_{j,n+1}^{(\pm)} x_{i, m}^{(\pm)}, \\
&&{\rm Sym}_{m_1,\cdots, m_{r}} \sum_{s=0}^{r}(-1)^s 
{[r]_i!\over [s]_i! [r-s]_i!}\times\n
&&x^{(\pm)}_{i, m_1}\cdots x^{(\pm)}_{i, m_s}
 x^{(\pm)}_{j,n}x^{(\pm)}_{i, m_{s+1}}\cdots x^{(\pm)}_{i, m_r}=0, 
\quad(i\ne j, r=1-a_{ij})
\end{eqnarray}
where  $\phi^\pm_{i,\mp r}$ $(r\ge 0)$  are expressed 
 in terms of $k^{\pm 1}_i$ and $h_{i,r}$'s by (\ref{eq:phi}).

\begin{thm}[\cite{Be1}] \label{thm:be}
For nontwisted affine Kac--Moody algebra $g$, there exists
 an algebra isomorphism $\psi : \D(g)\to U_q^\Sigma(g)$
 determined by
\begin{eqnarray*}
&&x_{i,k}^{(\pm)}\mapsto x_{i,k}^{(\pm)}, 
\quad h_{i,r}\mapsto h_{i,r},\quad
k_i^{\pm 1}\mapsto k_{i}^{\pm 1},
\quad  C^{\pm 1}\mapsto k_{\pm \delta},
\quad D^{\pm 1}\mapsto k_{\pm \Lambda_0/d_0}.
\end{eqnarray*}
\end{thm}

\subsection{}
In this subsection, we shall write down  the 
multiplicative formula of $\Theta$, following  the approach 
by Beck \cite{Be1},~\cite{Be2} and Lusztig \cite{Lt}.
Most of the necessary results have already been obtained by them.

Firstly we shall describe their results.
Fix  $x=\comega_{j_m}\cdots \comega_{j_1}\in \bQ^\vee\subset W$
such that $(x|\alpha_i)>0$ $(1\le i\le l)$.
Let $s_{i_n}\cdots s_{i_1}$ be its 
reduced expression obtained by the  concatenation of 
the reduced presentations  of $\comega_{j_k}$'s 
and canceling the elements of  $\Pi$.
Let  us define $\beta_k$ $(k\in \Z)$ by 
$$
\beta_{k}=s_{i_1}\cdots s_{i_{k-1}}(\alpha_{i_k}), \quad 
\beta_{k+rn}=x^{-r}(\beta_k)\, (r\ge 0),
\quad \beta_{k-rn}=x^{r}(-\beta_k) \, (r > 0)
$$
where $1\le k\le n$.  Then they are distinct 
and run over the whole set of   positive real roots of $g$.
Moreover $\beta_k$ is a   positive (resp.  negative) 
 root of $\bg$ mod $\Z\delta$ if $k>0$ (resp. $k\le 0$).
We define  $E_k\in U^+_{\beta_k}$ and   
$F_k\in U^-_{-\beta_k}$ $(k\in \Z)$ by
\begin{eqnarray}
E_{k+rn}&=&T_{x}^{-r}T_{i_1}^{-1}\cdots T_{i_{k-1}}^{-1}(e_{i_k})
\quad (1\le k\le n, r\ge 0)\n
E_{k-rn}&=&k_{\beta_{k-rn}}T_{x}^{r}T_{i_1}^{-1}\cdots
 T_{i_{k-1}}^{-1}(f_{i_k})\quad (1\le k\le n, r>0)\n
F_k&=&\Omega(E_k)
\end{eqnarray}
and set
\begin{eqnarray}
&&E^{{\bf m}}=(E_1^{m_1}E_2^{m_2}\cdots ) 
\prod_{1\le i\le l\atop r>0} h_{i,r}^{m_{i,r}}
(\cdots E_{-1}^{n_{-1}} E_{0}^{m_0})\n
&&F^{{\bf m}}=(F_1^{m_1}F_2^{m_2}\cdots ) 
\prod_{1\le i\le l\atop r>0} h_{i,-r}^{m_{i,r}}.
(\cdots F_{-1}^{m_{-1}} F_{0}^{m_0})
\end{eqnarray}
Here $h_{i, \pm r}\in U^\pm_{\pm r\delta}$
 $(r>0)$ is defined in (\ref{eq:phi}), and  
${\bf m}=(m_k, m_{j,r})_{k\in \Z, r>0, 1\le j\le l}$
 is a family   of  nonnegative integers such that 
all vanish except for a finite number of them.
For $k\in \Z$  we let $\bar{k}$ denote the integer  
which is between $1$ and $n$, and equal to $k$ mod $n$.
Finally, 
let  $(\,,\,): U^{\ge 0}\times U^{\le 0}\rightarrow F$
be the Hopf pairing \cite{Dr},~\cite[Proposition 2.1.1]{T}
 determined by  bilinearity and the following,
\begin{eqnarray}
&& (k_\lambda,k_\mu)=q^{-(\lambda|\mu)},
\quad (e_i,k_\lambda)=(k_\lambda,f_i)=0,\quad
(e_i,f_j)=-\delta_{ij}/(q_i-q_i^{-1}),\n
&&(x,y_1y_2)=(\Delta(x),y_1\otimes y_2),
\quad (x_1x_2,y)=(x_2\otimes x_1, \Delta(y))\n
&&(x, x_1, x_2\in U^{\ge 0},\quad y, y_1, y_2\in U^{\le 0})
\end{eqnarray}

\cite[Proposition 3]{Be2} and \cite[Proposition 40.2.4]{Lt}  give
\begin{pr}\label{pr:real}\cite{Be2}, \cite{Lt}

(1)$\{E^{{\bf m}}\}$ and $\{F^{{\bf m}}\}$
 form a basis of $U^+$ and $U^-$, respectively.

(2)$
(E^{{\bf m}}, F^{{\bf n}})=
(\prod_{{1\le i \le l\atop r>0}} h_{i,r}^{m_{i,r}},
\prod_{{1\le i\le l\atop r>0}} h_{i,-r}^{n_{i,r}})\times
\prod_{k\in \Z}\delta_{m_k, n_k}(e_{i_{\bar{k}}}^{m_k},
f_{i_{\bar{k}}}^{m_k}).
$
\end{pr}

Therefore we have only to show the following. 
For $r\in \Z_{>0}$, let 
$(C_{ij}(r))_{1\le i,j\le l}$  be  the inverse matrix of 
${\displaystyle \left({[ra_{ij}]_i(q-q^{-1})
\over r(q_j-q_j^{-1})}\right)_{1\le i,j \le l}}$
and set
$\th_{i,-r}=\sum_{j=1}^lC_{ij}(r)h_{j,-r}$.

\begin{lm}\label{lm:imaginary}
$$
(\prod_{i,r} h_{i,r}^{m_{i,r}},\prod_{i,r} \th_{i,-r}^{n_{i,r}})
=\prod_{i,r}\delta_{m_{i,r}, n_{i,r}} m_{i,r}!
\left({-1\over q-q^{-1}}\right)^{m_{i,r}}
$$
\end{lm}
Proof.
Let $\la\, ,\,\ra : U\times U\rightarrow F(q^{1/2})$
 be the Killing form introduced in \cite{T}.
This form is bilinear and has the following properties
\begin{eqnarray}
&&\la \ad(u)\cdot v_1,v_2\ra=\la v_1,v_2\cdot \tad(u)\ra 
\quad ( u,v_1,v_2\in U) \label{eq:ad}\\
&&\la x k_\lambda,y k_\mu\ra=(x,y)q^{-(\lambda|\mu)/2}
\quad (x\in U^+, y\in U^-)
\label{eq:lambda}\\
&&\la U^{\ge 0}U^-_{-\mu},U^{\le 0}\ra=\{0\}\quad (\mu>0)
\end{eqnarray}
In particular, (\ref{eq:lambda}) implies 
that $\la x, y k_\lambda\ra$
 is independent of $\lambda\in \Sigma$ 
if $x\in U^+$ and $y\in U^{\le 0}$.
Since $\th_{i,-r}\in U^-_{-r\delta}$ \cite{Be1},
$$
\Delta(\th_{i,-r})=1\otimes \th_{i,-r}+\th_{i,-r}
\otimes k_{-r\delta}\quad {\rm mod} \quad 
\mathop{\textstyle\bigoplus}_{r\delta>\beta>0}
U^-_{-(r\delta-\beta)}\otimes U^-_{-\beta}U^0
$$
Noting the above and (\ref{eq:h}),  we obtain from (\ref{eq:ad})
with $u=\th_{j,-s}$, $v_1=\prod_{i, r} h_{i,r}^{m_{i,r}}$ and 
$v_2=\prod_{i, r} \th_{i,-r}^{n_{i,r}}\cdot  k_\lambda$,

\begin{eqnarray*}
&&(q^{-(\lambda|s\delta)}-1)
\{
(\prod_{i,r} h_{i,r}^{m_{i,r}},\th_{j,-s}
\prod_{i,r} \th_{i,-r}^{n_{i,r}})\\
&&+{m_{j,s}\over q-q^{-1}}(h_{j,s}^{m_{j,s}-1}
\prod_{(i,r)\ne(j,s)} h_{i,r}^{m_{i,r}},
\prod_{i,r} \th_{i,-r}^{n_{i,r}})
\}
=\sum_{s\delta>\beta \ge 0}q^{-(\lambda|\beta)}x_\beta 
\end{eqnarray*}
Here  $x_\beta$'s are the elements of  $F$ independent of $\lambda$.
Since the above equality holds for any $\lambda\in \Sigma$,
the inside of the curly bracket  vanishes. 
From this,  the lemma follows.

\vskip2mm

The following is essentially due to Beck and Lusztig,
 cf. \cite{LSS},~\cite{KT1},~\cite{KT2},~\cite{KT3}.
Set
\begin{eqnarray}
&&\Theta^{0}=\exp\Bigl(-(q-q^{-1})
\sum_{1\le i,j\le l\atop r>0}C_{ij}(r) h_{i,r}
\otimes h_{j,-r}\Bigr),\n
&&\Theta^+=\Theta_1\Theta_2\cdots,
 \quad \Theta^-=\cdots\Theta_{-1}\Theta_0, 
\end{eqnarray}
where
$$
\Theta_k=\sum_{m\ge 0}(-1)^mq_i^{-{m(m-1)\over 2}}
{(q_i-q_i^{-1})^m\over [m]_i!} E_k^m\otimes F_k^m
\quad (i=i_{\bar{k}}).
$$

\begin{pr}\label{pr:mult}

(1) $\Theta=\Theta^+\Theta^0\Theta^-$.

(2)
Let $i\in \{1,\cdots,l\}$.
If  $\beta_k=\alpha_i+m\delta$  ($k>0$, $m\ge 0$),
 then $E_k\otimes F_k=x^{(+)}_{i,m}\otimes x^{(-)}_{i,-m}$.
If  $\beta_k=-\alpha_i+m\delta$  ( $k\le 0$, $m>0$),
 then $E_k\otimes F_k=\Psi(x^{(-)}_{i,m}\otimes x^{(+)}_{i,-m})$.
\end{pr}

Proof.
(1) Since $\Theta_\mu$  is the canonical  element 
of $(\,,\,)|_{U_\mu^+\times U_{-\mu}^-}$ \cite{Dr},~\cite{T},
the claim follows from Proposition \ref{pr:real} 
and Lemma \ref{lm:imaginary}.

(2) Let $x_i$ denote $e_i$ or $f_i$.
When  $1\le k\le n$,
 $T_{i_1}^{-1}\cdots T_{i_{k-1}}^{-1}(x_{i_k})
=T_{\comega_{j_1}}^{-1}\cdots T_{\comega_{j_{r-1}}}^{-1}$
 $\times T_{p_1}^{-1}\cdots T_{p_{s-1}}^{-1}(x_{p_s})$
 for some $r$, $s$, where
$\pi s_{p_t}\cdots s_{p_1}$ is the  reduced expression
 of $\comega_{j_r}$. Since the set
 $\{s_{p_1}\cdots s_{p_{s-1}}(\alpha_{p_s})|1\le s\le t\}$ contains 
only $\alpha_{p_1}=\alpha_{j_r}$ as a root  of the 
form $\pm \alpha_i+m\delta$ $(1\le i\le l, m\ge 0)$, (2) follows from 
$T_{\omega_i^\vee}T_{\omega_j^\vee}
=T_{\omega_j^\vee}T_{\omega_i^\vee}$ and 
$T_{\omega_i^\vee}x_j=x_j$ $(1\le i\ne j\le l)$ \cite{Be1}.


\section{$L$ operators and  Drinfeld's generators}

In the rest of the paper, 
we shall consider only $g=A_l^{(1)}, B_l^{(1)}, C_l^{(1)}, D_l^{(1)}$.
Let $L^\pm(z)=L^{\pm, u}(z)L^{\pm,d}(z)L^{\pm, l}(z)$
 be the unique decomposition of the $L$ operators
such that $L^{\pm,a}(z)\in (M_N(F)\otimes \A)[[z^{\pm 1}]]
\simeq M_N(\A[[z^{\pm 1}]])$ $(a=u,d,l)$,
$L^{\pm,d}(z)$  is diagonal, and
$L^{\pm, u}(z)$ (resp. $L^{\pm, l}(z)$) is  upper (resp. lower)
triangular with  identity elements on the diagonal  \cite{DF}.
For $\Theta^a$ $(a=\pm,0)$,  define 
$((\Theta^{a})^{\pm 1})_{\mu}$ $(\mu\in Q_+)$ as for $\Theta$. 
Let us introduce  
$\hL^{\pm,a}(z)\in (M_N(F)\otimes \U)[[z^{\pm 1}]]$ $(a=u,d,l)$ by
\begin{eqnarray}
&&\hL^{+,u}(z)=\hLL^{+,u}(z),\quad \hL^{+,d}(z)=\hLL^{+,d}(z)\hT^{-1},
\quad \hL^{+,l}(z)=\hT\hLL^{+,l}(z)\hT^{-1},\n
&&\hL^{-,u}(z)=\hT\hLL^{-,u}(z)\hT^{-1},\quad 
\hL^{-,d}(z)=\hT\hLL^{-,d}(z),\quad \hL^{-,l}(z)=\hLL^{-,l}(z),
\end{eqnarray}
where
$$
\hLL^{+,a}(z)=\sum_{\mu\ge 0} (\rho \otimes 1)
(\Theta^{a'}_\mu) z^{{(\Lambda_0|\mu)\over d_0}},\quad 
\hLL^{-,a}(z)=\sum_{\mu\ge 0} (\rho \otimes 1)\left(
((\Theta^{a''}_{21})^{-1})_\mu\right)
 z^{-{(\Lambda_0|\mu)\over d_0}},
$$
and $u'=l''=+$, $d'=d''=0$ and $l'=u''=-$ in the last equations.
Noting the property of $\beta_k$ as a root of 
$\bg$ mod $\Z\delta$ and the weight of $h_{i, r}$, we obtain 
\begin{lm}
If we identify $\A$ with $\U$ by the map 
$\phi$, $L^{\pm,a}(z)=\hL^{\pm,a}(z)$. 
\end{lm}

This, together with Proposition \ref{pr:mult}
 and Theorem \ref{thm:be}, gives the relation between 
the $L$ operators and Drinfeld's  generators.
 Before we state the results, we prepare some notations.
For $1\le i\le l$ , set 
\begin{eqnarray}
&&e_i^+(z)=\sum_{m>0}x^{(+)}_{i,-m}z^m,\quad 
f_i^+(z)=\sum_{m\ge 0}x^{(-)}_{i,-m}z^m,\quad 
\phi^+(z)=\sum_{m\ge 0}\phi^+_{-m}z^{m},\n
&&e_i^-(z)=\sum_{m\ge 0}x^{(+)}_{i,m}z^{-m}, \quad 
f_i^-(z)=\sum_{m>0}x^{(-)}_{i,m}z^{-m},\quad 
\phi^-(z)=\sum_{m\ge 0}\phi^-_{m}z^{- m}.\n
&&
\end{eqnarray}
Let $\rho_\pm(x_i)$ $(x=e,f)$ denote 
the expressions $x_i^\pm$ in Appendix A.2.
Let further $\epsilon(i)$ ($1\le i\le N$) denote $i$ 
for $\epsilon=+$ and $N+1-i$ for $\epsilon=-$.
For $L^{\pm,a}(z)$ $(a=u,l)$, 
let  $\equiv$ denote the equality modulo $\oplus'_{i,j} 
(F E_{ij}\otimes \A)[[z^{\pm 1}]]$, where the sum is taken over 
$i$, $j$ such that $\eta_i-\eta_j\ne 0, \epsilon  \alpha_k$
 ($1\le k\le l$, $\epsilon=+$ for $a=u$, $=-$ for $a=l$).
Finally, define $b\in F^\times $ by
$-o(l)b=q^{-(l+1)}$ for $g=A_l^{(1)}$, 
$=x/y\times q_0^{-h^\vee\!/2}$ (see (\ref{eq:xy}) for $x$ and $y$)
for $g=B_l^{(1)}$ and $=q_0^{-h^\vee\!/2}$ 
 for $g=C_l^{(1)}$, $D_l^{(1)}$.
($b$  is an inessential constant determined by 
the correspondence between  the $e_j$ and $x^{(-)}_{i,1}$.)

By direct calculations,  we obtain 
\begin{thm}
Let $g$ be $A_l^{(1)}$, $B_l^{(1)}$, $C_l^{(1)}$ or  $D_l^{(1)}$.
If we identify $\D(g)$ with the subalgebra of $\A(g)$ via 
 $(\tau_{b}\circ \phi)^{-1}|_{U_q^\Sigma(g)} \circ \psi$,
  the following holds.

(1) $g=A_l^{(1)}${\rm :}
\begin{eqnarray*}
\mp(L^{\pm,u}(z)-1)&\equiv&(q-q^{-1}) \sum_{1\le i\le l}
\rho(e_i)\otimes f_i^\pm(q^i C^{-{1\mp1\over 2}}z)\\
\mp(L^{\pm,l}(z)-1)&\equiv&(q-q^{-1}) \sum_{1\le i\le l}
\rho(f_i)\otimes e_i^\pm(q^i C^{-{1\pm 1\over 2}}z)\\
L_{ii}^{\pm ,d}(z)L_{i+1\, i+1}^{\pm ,d}(z)^{-1}&=
&\phi^\pm_i(q^iz)\quad (1\le i\le l)
\end{eqnarray*}

(2) $g=B_l^{(1)}, C_l^{(1)}, D_l^{(1)}${\rm :}
\begin{eqnarray*}
\mp(L^{\pm,u}(z)-1)&\equiv&
\sum_{1\le i\le l-1\atop \epsilon=\pm}(q_i-q_i^{-1})
 \rho_\epsilon(e_i)\otimes 
f_i^\pm(q_i^{\epsilon i} q_0^{-\epsilon h^\vee\!/2} 
C^{-{1\mp1\over 2}}z)\\
+ (q_l-q_l^{-1})&\times& \left\{
\begin{array}{rl}
\sum_{\epsilon=\pm}\rho_\epsilon(e_l)\otimes 
f_l^\pm(q_l^\epsilon C^{-{1\mp1\over 2}}z)&\quad 
\mbox{ for $g=B_l^{(1)}$}\\
\rho(e_l)\otimes f_l^\pm(C^{-{1\mp1\over 2}}z)&
\quad \mbox{ for $g=C_l^{(1)}, D_l^{(1)}$}
\end{array}\right.\\
\mp(L^{\pm,l}(z)-1)&\equiv&
\sum_{1\le i\le l-1\atop \epsilon=\pm}(q_i-q_i^{-1})
 \rho_\epsilon(f_i)\otimes 
e_i^\pm(q_i^{\epsilon i} q_0^{-\epsilon h^\vee\!/2} 
C^{-{1\pm 1\over 2}}z)\\
+ (q_l-q_l^{-1})&\times& \left\{
\begin{array}{rl}
\sum_{\epsilon=\pm}\rho_\epsilon(f_l)\otimes 
e_l^\pm(q_l^\epsilon C^{-{1\pm 1\over 2}}z)&\quad 
\mbox{ for $g=B_l^{(1)}$}\\
\rho(f_l)\otimes e_l^\pm(C^{-{1\pm 1\over 2}}z)&
\quad \mbox{ for $g=C_l^{(1)}, D_l^{(1)}$}
\end{array}\right.\\
L_{\epsilon(i)\, \epsilon(i)}^{\pm ,d}(z)^\epsilon
L_{\epsilon(i+1)\, \epsilon(i+1)}^{\pm ,d}(z)^{-\epsilon}
&=&\phi^\pm_i(q_i^{\epsilon i}q_0^{-\epsilon h^\vee\!/2}z)
\quad (1\le i\le l-1, \epsilon=\pm)\\
L_{\epsilon(l)\, \epsilon(l)}^{\pm ,d}(z)^\epsilon
L_{\epsilon(l+1)\, \epsilon(l+1)}^{\pm ,d}(z)^{-\epsilon}
&=&\phi^\pm_l(q_l^{\epsilon }z)\quad(\epsilon=\pm)\quad  
\mbox{for $g=B_l^{(1)}$}\\
L_{l l}^{\pm ,d}(z)L_{l+1\, l+1}^{\pm ,d}(z)^{-1}
&=&\phi^\pm_l(z)\quad \mbox{for $g=C_l^{(1)}$}\\
L_{l-1\, l-1}^{\pm ,d}(z)L_{l+1\, l+1}^{\pm ,d}(z)^{-1}
&=&
L_{l l}^{\pm ,d}(z)L_{l+2\, l+2}^{\pm ,d}(z)^{-1}=
\phi^\pm_l(z)\quad \mbox{for $g=D_l^{(1)}$}
\end{eqnarray*}

\end{thm}

\vskip5mm

\noindent {\it Remark.} In \cite{KT3},  another Hopf 
structure of quantum group found by  Drinfeld is shown to be
obtained by the twisting by $\F=(\Theta^+_{21})^{-1}$
 (in our notation). The above correspondence 
between the  $L$ operators and Drinfeld's generators
gives an intuitive  explanation  
of its comultiplication formula  as follows.
After the twist, $\Theta^\F=\Theta^0\Theta^-\Psi(\Theta_{21}^+)$. 
Define $\hL^{\pm,\F}(z)$ for $\Theta^\F$
as before. Then  $\hL^{+,\F}(z)=\hL^{+,d}(z)\hL^{+,l}(z)
\hL^{-,l}(k_{-\delta}z)^{-1}$ and 
$\hL^{-,\F}(z)=\hL^{+,u}(k_{-\delta}z)^{-1}\hL^{-,u}(z)\hL^{-,d}(z)$.
These  have the lower and upper triangular forms, respectively,
and satisfy the comultiplication formula in (\ref{eq:copr})
for $\Delta^\F(\cdot)=\F\Delta(\cdot)\F^{-1}$.  From this,
the comultiplication formula for Drinfeld's generators easily follows.
The above argument is formal, but the justification would  be possible.


\appendix
\section{Appendix   $\quad$ $U'$ module  $(\rho, V)$  and $w$}

\subsection{ $N$ and  $\eta=(\eta_1,\cdots, \eta_N)$}

We set $\epsilon_i=p\omega_1-(\alpha_1+\cdots\alpha_{i-1})$
 $(1\le i\le l)$.
(Note that in the case $g=D_l^{(1)}$, $v_l$ and $v_{l+1}$
 can be interchanged.)

\begin{eqnarray*}
&&(1)\,A_l^{(1)}\,(l\ge 1): N=l+1, \quad 
\eta=(\epsilon_1,\cdots,\epsilon_l, -\sum_{i=1}^l\epsilon_i)\\
&&(2)\, B_l^{(1)}\,(l\ge 3), A_{2l}^{(2)}\, (l\ge 1):
N=2l+1, \quad \eta
=(\epsilon_1,\cdots,\epsilon_l,0,-\epsilon_l,\cdots,-\epsilon_1)\\
&&(3)\, C_l^{(1)}\,(l\ge 2), A_{2l-1}^{(2)}\,(l\ge 3), 
D_l^{(1)}\,(l\ge 4): N=2l, \quad \eta
=(\epsilon_1,\cdots,\epsilon_l,-\epsilon_l,\cdots,-\epsilon_1)\\
&&(4)\, D_{l+1}^{(2)}\,(l\ge 2):N=2l+2, \quad 
\eta=(\epsilon_1,\cdots,\epsilon_l,0,-\epsilon_l,\cdots,-\epsilon_1,0)
\end{eqnarray*}

\subsection{ $(\rho, V)$}
We give the matrix representations of $e_i$ and $f_i$
 $(0\le i\le l)$ with repect to the basis $v_1,\cdots, v_N$. 
$E_{ij}$ is a matrix unit and $\t f_i$ denotes 
the transpose of $f_i$. $x_i^{B_l}$  and $x_i^{C_l}$ 
($x=e,f$, $1\le i\le l$) below stand for the expressions
 $x_i$ given for $B_l^{(1)}$ and $C_l^{(1)}$, respectively.

\begin{eqnarray}
(1)\, A_l^{(1)}:&&  e_i=E_{i\,i+1}=\t f_i
\quad (1\le i\le l),\quad  e_0=E_{N1}=\t f_0\\
(2)\, B_l^{(1)},&& \hskip-6mm C_l^{(1)}, D_l^{(1)}:\n
({\rm i})\, 1\le i <l:&&e_i=e_i^{+}+e_i^{-},
\quad f_i=f_i^{+}+f_i^{-},\n
&&e_i^{+}=E_{i\,i+1}=\t f_i^{+},
\quad e_i^{-}=-E_{N-i\,\,N+1-i}=\t f_i^{-},\\
({\rm ii})\,  i=l,0:&&\n
B^{(1)}_l:&&e_l=e_l^{+}+e_l^{-},
\quad f_l=f_l^{+}+f_l^{-},\n
&&x^{-1}e_l^{+}=E_{l\,l+1}=y^{-1} \t f_l^{+},
\quad x^{-1} e_l^{-}=-E_{l+1\,l+2}=y^{-1}\t f_l^{-}\n
&&e_0=E_{N-1\, 1}-E_{N\,2}=\t f_0\n
&&x,y \in A \mbox{  such that  } xy=[2]_l\label{eq:xy}\\
C^{(1)}_l:&&e_l=E_{l\, l+1}=\t f_l,\quad e_0=E_{N 1}=\t f_0\\
D^{(1)}_l:&&e_l= E_{l-1\, l+1}-E_{l\,l+2}=\t f_l,\n
&&e_0=E_{N-1\, 1}-E_{N 2}=\t f_0\\
(3)\, A_{2l}^{(2)}:&& x_i=Mx^{B_l}_iM^{-1}
\quad (x=e,f,\quad 1\le i\le l)\n
&&e_0=E_{N 1}=\t f_0\n
&&M=\diag(\overbrace{1\cdots 1}^{l+1},-1,1,\cdots, (-1)^l)\\
(4)\, A_{2l-1}^{(2)}:&& x_i=Mx^{C_l}_iM^{-1}
\quad (x=e,f, \quad 1\le i\le l)\n
&&e_0=E_{N-1\, 1}+E_{N 2}=\t f_0\n
&&M=
\diag(\overbrace{1\cdots 1}^{l+1},-1,1,\cdots, (-1)^{l-1})\\
(5)\, D_{l+1}^{(2)}:&& x_i=x^{B_l}_i
\quad  (x=e,f, \quad 1\le i\le l)\n
&&x^{-1}e_0=E_{N 1}+E_{N-1\, N}=y^{-1}\t f_0\n
&&x,y\in A \mbox{  such that  } xy=[2]_0
\end{eqnarray}

\subsection{ $w$}

\begin{eqnarray*}
(1)\,  g=A_l^{(1)}:&&
w=\sum_{\gamma\in {\cal S}_{l+1}}
 (-q)^{l(\gamma)}v_{\gamma(1)}\otimes\cdots\otimes 
v_{\gamma(l+1)}.\\
&&(l(\gamma): \mbox{length function})\n
(2)\,  g\ne A_l^{(1)}:&&
 w=\sum_{1\le i,j\le N}J_{ij}v_i\otimes v_j,\n
J_{ij}&=&
u_i\delta_{i+j, N+1} \mbox{  for  $g\ne D_{l+1}^{(2)}$},
=u_i\delta_{i+j, N}+u_N\delta_{i N}\delta_{jN}
\mbox{   for  $g= D_{l+1}^{(2)}$}.\n
\mbox{List of}&& \hskip-7mm(u_1,\cdots u_N)\n
({\rm i})\, B_l^{(1)}:&&\hskip-3mm
 (q^{-(2l-1)},\cdots, q^{-3},q^{-1},q,q,q^3,\cdots, q^{2l-1}) \\
({\rm ii})\, C_l^{(1)}:
&&\hskip-3mm(-q^{-l},\cdots,-q^{-2}, -q^{-1},q,q^2,\cdots, q^l) \\
({\rm iii})\, D_l^{(1)}:
&&\hskip-3mm (q^{-(l-1)},\cdots,q^{-1},1,1,q, \cdots, q^{l-1})\\
({\rm iv})\, A_{2l}^{(2)}:&&\hskip-3mm
 ((-1)^lq^{-(2l-1)},\cdots,q^{-3},-q^{-1},q,-q,q^3,
\cdots, (-1)^l q^{2l-1})\\
({\rm v})\, A_{2l-1}^{(2)}:&&\hskip-3mm
((-1)^lq^{-l},\cdots,q^{-2},-q^{-1},q,-q^2,\cdots, (-1)^{l-1}q^{l})\\
({\rm vi})\, D_{l+1}^{(2)}:&&\hskip-3mm 
(q^{-(2l-1)},\cdots,q^{-3},q^{-1},q,q,q^3, \cdots, q^{2l-1},-q)
\end{eqnarray*}


\end{document}